\documentclass[11pt]{amsart}

\usepackage{amssymb,amsthm,amsmath}
\usepackage[numbers,sort&compress]{natbib}
\usepackage{color}
\usepackage{graphicx}
\usepackage{tikz}
\usepackage[colorlinks,
            linkcolor=blue,
            anchorcolor=blue,
            citecolor=blue
            ]{hyperref}

\let\newpf\proof \let\proof\relax 
\newenvironment{pf}{\newpf[\proofname]}{\qed\endtrivlist}

\newcommand{\ba}{\overline{A}}

\def\be{\begin{equation}}
\def\ee{\end{equation}}

\def\ba{{\begin{align}}}
\def\ea{{\end{align}}}

\def\bm{\begin{matrix}}
\def\em{\end{matrix}}

\def\0{{\mathbf 0}}

\newtheorem{Theorem}{Theorem}[section]
\newtheorem*{Theorem*}{Theorem 1.1}
\newtheorem{Lemma}{Lemma}[section]
\newtheorem{Proposition}{Proposition}[section]

\newtheorem{Remark}{Remark}[section]

\numberwithin{equation}{section}

\theoremstyle{definition}

\newcommand{\C}{{\mathbb C}}

\newcommand{\N}{{\mathbb N}}

\newcommand{\R}{{\mathbb R}}
\newcommand{\T}{{\mathbb T}}

\newcommand{\Z}{{\mathbb Z}}

\newcommand{\la}{\langle}
\newcommand{\ra}{\rangle}

\def\B0{{\bold{0}}}

\catcode`\@=12

\def\Empty{}
\newcommand\oplabel[1]{
  \def\OpArg{#1} \ifx \OpArg\Empty {} \else
    \label{#1}
  \fi}

\newcommand{\comm}[1]{}
\newcommand{\comment}[1]{}

\begin{document}

\title[Sharp H\"older continuity  of  Lyapunov exponent]{ Sharp H\"older continuity of the Lyapunov exponent of finitely differentiable quasi-periodic cocycles  }

\author{Ao Cai} \address{
Department of Mathematics, Nanjing University, Nanjing 210093, China
}

 \email{godcaiao@126.com}

\author{Claire Chavaudret} \address{
Laboratoire J.A. Dieudonn\' e, Universit\' e de Nice-Sophia Antipolis (Parc Valrose), 06108 Nice Cedex 02, France
}

 \email{claire.chavaudret@unice.fr}

\author {Jiangong You}
\address{
Chern Institute of Mathematics and LPMC, Nankai University, Tianjin 300071, China} \email{jyou@nankai.edu.cn}

\author{Qi Zhou}
\address{
Department of Mathematics, Nanjing University, Nanjing 210093, China
}

 \email{qizhou628@gmail.com, qizhou@nju.edu.cn}

\setcounter{tocdepth}{1}

\begin{abstract}
We show that  if the base frequency is Diophantine, then the Lyapunov exponent of a $C^{k}$ quasi-periodic $SL(2,\R)$ cocycle is  $1/2$-H\"older continuous  in the almost reducible regime, if $k$ is large enough.  As a consequence,  we show  that if the  frequency is Diophantine, $k$ is large enough,  and the  potential is  $C^k$ small, then   the integrated density of states of  the  corresponding quasi-periodic Schr\"{o}dinger operator  is  $1/2$-H\"older continuous.
\end{abstract}

\maketitle

\section{Introduction}

Let $(M,\mathcal{B},\mu)$ be a probability space and $f: M\rightarrow M$ be an invertible  measure preserving map, assuming $\mu$ is ergodic. Given a
measurable function  $A: M \rightarrow
SL(2,\R)$, the  linear cocycle given by $A$ over the base dynamics $f$ is the transformation:
$$
(f,A):M \times\R^{2} \rightarrow M \times\R^{2};(\theta,v)\mapsto (f(\theta),A(\theta)\cdot v).
$$
The iterates of $(f,A)$ have the form $(f,A)^n=(f^n,A_n)$, where
$
A_n(\theta)= A(f^{n-1}(\theta))\cdots A(f(\theta))A(\theta), n\geqslant 1 $ and $A_{-n}(\theta)= A_n(f^{-n}(\theta))^{-1}.$
The Lyapunov exponent is given by the formula
$$
L(f,A)=\lim_{n \to +\infty}\frac{1}{n}\int_{M} \ln \| A_n(\theta)\| d\mu.
$$

Lyapunov exponents appear naturally in the study of smooth dynamics. Continuity of Lyapunov exponents depends sensitively on  the finer topologies of the fiber and
 the base dynamics. This
 has been the object of considerable recent interests, see Viana \cite{V2} and references therein.  It is well known that in $C^0$ topology, discontinuity of Lyapunov exponent occurs at every non-uniformly hyperbolic cocycle, see \cite{furman,knill,thouvenot}. Moreover, Bochi \cite{Bochi1,Bochi2} proved that with an ergodic base system, any non-uniformly hyperbolic $SL(2,\R)$-cocycle can be approximated by cocycles with zero Lyapunov exponents in the $C^0$ topology. Therefore, if one wants to obtain some regularity of the Lyapunov exponent, finer topology than $C^0$ topology is necessary.

Now we discuss the impact of the base dynamics.  If the base dynamics has some hyperbolicity, then the Lyapunov exponent is continuous. For example,  Bocker and Viana \cite{BocV} proved continuity of Lyapunov exponents  for random products of $SL(2,\R)$ matrices  in the Bernoulli setting.  In higher dimensions, continuous dependence of all Lyapunov exponents for i.i.d. random products of matrices in $GL(d, R)$ was proved by Avila, Eskin, and Viana \cite{AEV}. If the base dynamics is a subshift of finite type or, more generally, a hyperbolic set, then Backes, Brown, Butler \cite{BBB} proved that  Lyapunov exponents are always continuous among H\"older continuous fiber-bunched $SL(2,\R)$-cocycles, which answers a famous conjecture  of Viana \cite{V2}.

 If the base is elliptic, things become more complicated: it will depend on the smoothness of $A(\theta)$ in a very sensitive way. Now assume that the base dynamics is quasiperiodic, i.e. $f: \T^d\rightarrow \T^d$ with $f(\theta)=\theta+\alpha$, where $\alpha$ is rationally independent.  If $A(\theta)$ is analytic, the Lyapunov exponent is always continuous \cite{B,BJ,JKS},  even in high dimensional $GL(d,\C)$ cocycles (see \cite{AJS}). However, if $A(\theta)$ is finitely differentiable,  Wang and You \cite{wangyou1} constructed  examples of discontinuity of the Lyapunov exponent in smooth quasi-periodic $SL(2,\R)$ cocycles even when $\alpha$ is of bounded type.

In this paper, we will show that if $A\in C^{k}(\T^{d},SL(2,\R))$, Lyapunov exponent can still be continuous, even H\"older continuous provided that the cocycle is almost reducible. Before explaining our result precisely, we introduce some basic concepts.  Recall that $\alpha \in\R^d$ is called {\it Diophantine} if there are $\kappa>0$ and $\tau>d-1$ such that $\alpha \in {\rm DC}(\kappa,\tau)$, where
\begin{equation}\label{dio}
{\rm DC}(\kappa,\tau):=\left\{\omega \in\R^d:  \inf_{j \in \Z}\left|\la n,\omega  \ra - j \right|
> \frac{\kappa}{|n|^{\tau}},\quad \forall \, n\in\Z^d\backslash\{0\} \right\}.
\end{equation}
We also set ${\rm DC}:=\bigcup_{\kappa>0,\tau>d-1} {\rm DC}(\kappa,\tau)$.  Given two cocycles $(\alpha,A_1)$, $(\alpha,A_2)\in \T^d    \times C^{k'}(\T^d,SL(2,\R))$,  one says that they are $C^k$ conjugated if there exists $Z\in C^{k}(2\T^d, $ $SL(2,\R))$, such that $$
Z(\theta+\alpha)^{-1}A_1(\theta)Z(\theta)=A_2(\theta).
$$
We say $(\alpha,A)$ is  $C^{k^{'},k}$ almost reducible, if   $A\in C^{k'}(\T^d,SL(2,\R))$, and closure of its $C^k$ conjugacies contains a constant. Within the above concepts, our main result is the following:

\begin{Theorem}\label{thm1}
Let $\alpha \in DC(\kappa,\tau)$, there exists a numerical constant $D_0>0$, such that  if $(\alpha,A)$ is $C^{k^{'},k}$ almost reducible with $k^{'}>k\geqslant D_0\tau$, then for any continuous map $B:\T^{d} \rightarrow SL(2,\C)$, we have
$$\lvert L(\alpha,A)-L(\alpha,B)\rvert \leqslant C_0\lVert B-A\rVert_0^{\frac{1}{2}},$$
where $C_0$ is a constant depending on $d,\kappa,\tau$.
\end{Theorem}

Typical example of $SL(2,\R)$ cocycles comes from the  Schr\"{o}dinger cocycles, where
$$
A(\theta)=S_E^{\lambda V}=\begin{pmatrix}\lambda  V(\theta)-E & -1 \\ 1 & 0 \end{pmatrix}.
$$
They come from one dimensional quasi-periodic Schr\"{o}dinger operators on $\ell^{2}(\Z)$:
$$
(H_{\lambda V,\alpha,\theta}u)_n=u_{n+1}+u_{n-1}+\lambda V(\theta+n\alpha)u_{n}=Eu_{n}.
$$
since a solution $H_{\lambda V,\alpha,\theta}u=Eu$ satisfies $A(\theta+n\alpha)\begin{pmatrix} u_{n} \\ u_{n-1} \end{pmatrix}=\begin{pmatrix} u_{n+1} \\ u_{n} \end{pmatrix}$.

For  Schr\"{o}dinger operators $H_{\lambda V,\alpha,\theta}$, another important concept related to the Lyapunov exponent is the integrated density of states (IDS),  which is the function $N:\R \rightarrow [0,1]$ defined by
$$
N_{\lambda V,\alpha}(E)=\int_{\T^{d}}\mu_{\lambda V,\alpha,\theta}(-\infty,E]d\theta,
$$
where $\mu_{\lambda V,\alpha,\theta}$ is the spectral measure of $H_{\lambda V,\alpha,\theta}$. As a local version of our main result, we are able to show the following:

\begin{Theorem}\label{cor1}
Let $\alpha \in DC(\kappa,\tau)$, $V\in C^{k^{'}}(\T^{d},\R)$ with $k^{'}\geqslant D_0\tau$, where $D_0$ is a numerical constant.  Then there exist $\lambda_{0},k$ depending  on $V,d,\kappa,\tau,k^{'}$ such that if $\lambda<\lambda_0$, then we have the following:
\begin{enumerate}
\item For any $E\in\R$, $(\alpha,S_E^{\lambda V})$ is $C^{k^{'},k}$ almost reducible.
\item  $N_{\lambda V,\alpha}$ is $1/2$-H\"older continuous.
\end{enumerate}
\end{Theorem}

\begin{Remark}
We provide the first positive result of H\"{o}lder continuity of the Lyapunov exponent for  finitely differentiable $SL(2,\R)$ cocycles. This kind of  $1/2$-H\"{o}lder continuity  is sharp due to  the presence of spectral gaps.
\end{Remark}

\begin{Remark}
Our choice of $k^{'}$ and $k$ is not optimal, readers are invited to consult section \ref{3.5} for more discussions.
\end{Remark}

Before giving the rough ideas of the proof,  let us first  go over more history on the study of high order regularity of Lyapunov exponents. Just note that by the Thouless formula, (weak) H\"older continuity of the Lyapunov exponent is equivalent to (weak) H\" older continuity of the integrated density of states.

If the cocycle is analytic, one should distinguish zero Lyapunov exponent and positive Lyapunov exponent. In the positive Lyapunov exponent regime,  it starts with the work of Goldstein and Schlag \cite{goldsteinschlag1} where they developed some sharp version of large deviation theorems  for real analytic potentials with strong Diophantine frequency; they further developed Avalanche Principle and proved that  Lyapunov exponent is H\"{o}lder continuous (one-frequency) or weak H\"{o}lder continuous (multi-frequency). For the Almost Mathieu operator $H_{\lambda}=\lambda \cos2\pi(k\omega+\theta)+\triangle$, Bourgain \cite{bourgain} proved that for large $\lambda$ depending on Diophantine $\omega$, the Lyapunov exponent is $\frac{1}{2}-\epsilon$-H\"{o}lder continuous for any $\epsilon>0$.  Later, Goldstein and Schlag \cite{goldsteinschlag2} generalized Bourgain's result \cite{bourgain}, and proved that  if the potential  is in a small $L^{\infty}$ neighborhood of a trigonometric polynomial of degree $k$, then the IDS is H\"older $\frac{1}{2k}-\epsilon$-continuous for all $\epsilon>0$.

On the other hand, in the zero Lyapunov exponent regime, based on Eliasson's perturbative KAM scheme \cite{eli92}, Amor \cite{amor} obtained $1/2$-H\"{o}lder continuity of the IDS for quasi-periodic Schr\"odinger operator  with a Diophantine frequency and small potential.  Avila and Jitomirskaya \cite{avilajitomirskaya} further proved that if the cocycle is $C^\omega$ almost reducible (consult section \ref{asd} for definitions), then the Lyapunov exponent is $1/2$ H\"older continous. Thus our result generalized Avila and Jitomirskaya's result \cite{avilajitomirskaya} to the differentiable case.

 Notice that all the results stated above require Diophantine or strong Diophantine conditions. In fact, for small potential and generic frequencies, it is possible to show that the Lyapunov exponent is not  H\"{o}lder continuous \cite{avila}. A recent breakthrough belongs to Avila \cite{avila}: for one-frequency Schr\"{o}dinger operators with general analytic potentials and irrational frequency, Avila has built the fantastic global theory saying that the Lyapunov exponent is a $C^\omega$-stratified function of the energy.

In the lower regularity case, Klein \cite{klein} proved that for Schr\"{o}dinger operators with potentials in a Gevrey class, the Lyapunov exponent  is   weak H\"{o}lder continuous on any compact interval of the energy provided that the coupling constant is large enough, the frequency is Diophantine and the potential satisfies some transversality condition.
Recently, Wang and Zhang \cite{wangzhang} obtained the weak H\"{o}lder continuity of the Lyapunov exponent as function of energies, for a class of $C^2$ quasi-periodic potentials and for any Diophantine frequency.

For other related results, Avila and Krikorian \cite{avilakrikorian} have studied so-called monotonic cocycles which are a class of smooth or analytic cocycles non-homotopic to constant. They showed that the Lyapunov exponent is smooth or even analytic, respectively. Besides, for a 1-periodic function satisfying a Lipschitz monotonicity condition, Jitomirskaya and Kachkovskiy \cite{jitomirskayakachkovskiy} showed that the Lyapunov exponent is almost Lipschitz continuous.
\\

Finally, we will comment on the method of the proof.  Our method is purely dynamical, and the philosophy is that nice quantitative almost reducibility implies nice spectral applications.
This approach, which was first developed in \cite{eli92},  has been proved to be very fruitful \cite{A1, avilajitomirskaya,AYZ1,LYZZ}. To our purposes, the key is to have a nice control of the growth of the cocycles, i.e. to have a nice control of the $C^0$ norm of the conjugacy and the perturbation in the almost reducible scheme. Since our system is $C^k$ instead of analytic, we will first perform the KAM scheme in the analytic topology, then obtain the corresponding $C^k$ estimates by analytic approximations \cite{zehnder}, which is very classical in KAM theory.  In order to  do the analytic  approximation,   we need the strong almost reducibility results  (i.e. the cocycle is almost reducible in a fixed band, and this band is arbitrarily close to the initial band). Although this kind of result was  first obtained by Chavaudret \cite{chavaudret2}, however, her estimates there are too rough to have good spectral applications (readers can consult section \ref{3.5} for more discussions). Instead, we will perform a new KAM scheme which was developed by Leguil-You-Zhao-Zhou \cite{LYZZ}, where they initially developed the scheme to obtain the exponential asymptotics on the size of spectral gaps for almost Mathieu operators.  With the KAM scheme \cite{LYZZ} in hand, the key is to choose the suitable analytic approximation, to have sufficiently subtle quantitative estimates.

\section{Preliminaries}

 For a bounded
analytic (possibly matrix-valued) function $F$ defined on $ \{ \theta |  | \Im \theta |< h \}$, let
$
\|F\|_h=  \sup_{ | \Im \theta |< h } \| F(\theta)\| $ and denote by $C^\omega_{h}(\T^d,*)$ the
set of all these $*$-valued functions ($*$ will usually denote $\R$, $sl(2,\R)$
$SL(2,\R)$). Also we denote $C^\omega(\T^d,*)=\cup_{h>0}C^\omega_{h}(\T^d,*)$, and set  $C^{k}(\T^{d},*)$  to be  the space of $k$ times differentiable matrix-valued functions. The norms are defined as $$\lVert F \rVert _{k}=\sup_{\substack{
                             k^{'}\leqslant k,
                             \theta \in \T^{d}
                          }}\lVert \partial^{k^{'}}F(\theta) \rVert,
$$
and
$$
\lVert F \rVert _{0}=\sup_{\theta \in \T^{d}}\lVert F(\theta) \rVert.
$$
For any $N>0$, define the truncating operators $\mathcal{T}_N$ on $C^\omega(\T^d,*)$ as
$$
(\mathcal{T}_N{f})(\theta)=\sum_{k\in \Z^{d},\lvert k \rvert<N}\hat{f}(k)e^{i<k,\theta>}
$$
and $\mathcal{R}_N$ as
$$
(\mathcal{R}_N{f})(\theta)=\sum_{k\in \Z^{d},\lvert k \rvert\geqslant N}\hat{f}(k)e^{i<k,\theta>}.
$$
\subsection{Hyperbolicity and (almost) reducibility}\label{asd}
We say the cocycle $(\alpha, A)$ is $uniformly$ $hyperbolic$ if for every $\theta \in \T^d$, there exists a continuous splitting $\C^2=E^s(\theta)\oplus E^u(\theta)$ such that for some constants $C>0,c>0$, and for every $n\geqslant 0$,
$$
\begin{aligned}
\lvert A_n(\theta)v\rvert \leqslant Ce^{-cn}\lvert v\rvert, & v\in E^s(\theta),\\
\lvert A_n(\theta)^{-1}v\rvert \leqslant Ce^{-cn}\lvert v\rvert, & v\in E^u(\theta+n\alpha).
\end{aligned}
$$
This splitting is invariant by the dynamics, which means that for every $\theta \in \T^d$, $A(\theta)E^{\ast}(\theta)=E^{\ast}(\theta+\alpha)$, for $\ast=s,u$. In this case, it is clear that we have $L(\alpha,A)>0$.

Given two analytic cocycles $(\alpha,A_1)$, $(\alpha,A_2)\in \T^d \times C^{\omega}(\T^d,SL(2,\R))$,  they are analytically conjugated if there exists $Z\in C^{\omega}(2\T^d,SL(2,\R))$, such that
$$
Z(\theta+\alpha)^{-1}A_1(\theta)Z(\theta)=A_2(\theta).
$$
We call the cocycle $(\alpha,A)$ is $C^{\omega}$ almost reducible if the closure of its analytic conjugacy class contains a constant.

\subsection{Analytic approximation}
Assume $f \in C^{k}(\T^{d},sl(2,\R))$. By \cite{zehnder}, there exists a sequence $(f_{j})_{j\geqslant 1}$, $f_{j}\in C_{\frac{1}{j}}^{\omega}(\T^{d},sl(2,\R))$ and a universal constant $C^{'}$, such that
\begin{eqnarray} \nonumber \lVert f_{j}-f \rVert_{k} &\rightarrow& 0 , \qquad  j \rightarrow +\infty, \\
\label{2.1}\lvert f_{j}\rvert_{\frac{1}{j}} &\leqslant& C^{'}\lVert f \rVert_{k}, \\   \nonumber \lvert f_{j+1}-f_{j} \rvert_{\frac{1}{j+1}} &\leqslant& C^{'}(\frac{1}{j})^k\lVert f \rVert_{k}.
\end{eqnarray}
Moreover, if $k\leqslant k^{'}$ and $f\in C^{k^{'}}$, then properties (2.1) hold with $k^{'}$ instead of $k$. That means this sequence is obtained from $f$ regardless of its regularity (since $f_{j}$ is the convolution of $F$ with a map which does not depend on $k$).

\section{Quantitative estimates of almost reducibility}

The concept of $C^{k^{'},k}$ almost reducibility is merely qualitative and does not involve any quantitative estimates. However, the surprising thing is that one can really derive quite good quantitative  estimates, since this concept captures the  essential of applicability of local theories.   In order to give  almost reducibility theorem for $C^{k}$ quasi-periodic $SL(2,\R)$ cocycle,  we are not going to perform a KAM scheme straightly on a $C^k$ cocycle. Instead, we will deal with almost reducibility for analytic cocycles and turn it into a finitely differentiable analogue by virtue of analytic approximation. So we first treat the analytic case before dealing with the differentiable case.

\subsection{Decomposition along resonances}

In the following subsection, parameters $\rho,\epsilon,N,\sigma$ will be fixed; one will refer to the situation where there exists $n_\ast$ with $0<\lvert n_\ast\rvert \leqslant N$ such that
$$
\lvert 2\rho- <n_\ast,\alpha> \rvert< \epsilon^{\sigma},
$$
as the "resonant case". The integer vector $n_\ast$ will be referred to as a "resonant site".

\bigskip
Resonances are linked to a useful decomposition of the space $C^{\omega}_{r}(\T^{d},su(1,1))$ which is defined as follows:
for any given $\eta>0$, $\alpha\in \R^{d}$ and $A\in SU(1,1)$, we decompose $\mathcal{B}_r=C^{\omega}_{r}(\T^{d},su(1,1))=\mathcal{B}_r^{nre}(\eta) \bigoplus\mathcal{B}_r^{re}(\eta)$ in such a way that for any $Y\in\mathcal{B}_r^{nre}(\eta)$,
\begin{equation}\label{space}
A^{-1}Y(\theta+\alpha)A\in\mathcal{B}_r^{nre}(\eta), \qquad \lvert A^{-1}Y(\theta+\alpha)A-Y(\theta)\rvert_r\geqslant\eta\lvert Y(\theta)\rvert_r.
\end{equation}
Moreover, let $\mathbb{P}_{nre}$ and $\mathbb{P}_{re}$ be the standard projections from $\mathcal{B}_r$ onto $\mathcal{B}_r^{nre}(\eta)$ and $\mathcal{B}_r^{re}(\eta)$ respectively.

\bigskip
Within the above notations, the basic fact for our purpose is the following result:

\begin{Lemma}\label{lem2} Assume that $\epsilon\leqslant (4\lVert A\rVert)^{-4}$ and  $\eta \geqslant 13\lVert A\rVert^2{\epsilon}^{\frac{1}{2}}$. For any $g\in \mathcal{B}_r$ with $|g|_r \leqslant \epsilon$,  there exist $Y\in \mathcal{B}_r$ and $g^{re}\in \mathcal{B}_r^{re}(\eta)$ such that
$$
e^{Y(\theta+\alpha)}(Ae^{g(\theta)})e^{-Y(\theta)}=Ae^{g^{re}(\theta)},
$$
with $\lvert Y \rvert_r\leqslant \epsilon^{\frac{1}{2}}$ and $\lvert g^{re}\rvert_r\leqslant 2\epsilon$.
\end{Lemma}

\begin{Remark}\label{rem2}
The continuous version  of this result appears  in Lemma 3.1 of \cite{houyou}.  We will give a different proof using the quantitative Implicit Function Theorem in the Appendix.
Just point out that  the proof only relies on the fact that $\mathcal{B}_r$ is a Banach space, thus it also applies to $C^k$ and $C^0$ topology. \end{Remark}

\subsection{Analytic case} In this subsection, we concentrate on the following analytic quasi-periodic $SL(2,\R)$ cocycle:
$$
(\alpha,Ae^{f(\theta)}):\T^{d}\times\R^{2} \rightarrow \T^{d}\times\R^{2};(\theta,v)\mapsto (\theta+\alpha,Ae^{f(\theta)}\cdot v),
$$
where
$f\in C^{\omega}_{r_0}(\T^{d},sl(2,\R)), \ \ r_0>0, d\in \N^+$£¬
and $\alpha\in DC(\kappa,\tau)$.
Notice that $A$ has eigenvalues $\{e^{i\rho},e^{-i\rho}\}$ with $\rho \in \R\cup i\R$.

Through the proof, we find that analytic approximation requires us to show a result concerning a stronger sense of almost reducibility, i.e. we may shrink the analytic radius as little as we want in each iteration step and the perturbation still tends to zero.
Moreover, we must have a good control of the analytic norm of the conjugation map.

\begin{Proposition}\label{prop1}
Let $\alpha\in DC(\kappa,\tau)$, $\kappa,r>0$, $\tau>d-1$, $\sigma=\frac{1}{10}$.
Suppose that $A\in SL(2,\R)$, $f\in C^{\omega}_{r}(\T^{d},sl(2,\R))$.  Then for any $r'\in (0,r)$, there exist $c=c(\kappa,\tau,d)$ and a numerical constant $D$ such that if
\begin{equation}\label{estf}
\lvert f \rvert_r\leqslant\epsilon \leqslant \frac{c}{\lVert A\rVert^D}(r-r')^{D\tau},
\end{equation}
then there exist $B\in C^{\omega}_{r'}(2\T^{d},SL(2,\R))$, $A_{+}\in SL(2,\R)$ and $f_{+}\in C^{\omega}_{r'}(\T^{d},$
$sl(2,\R))$ such that
$$
B(\theta+\alpha)(Ae^{f(\theta)})B^{-1}(\theta)=A_{+}e^{f_+(\theta)}.
$$
More precisely, let $N=\frac{2}{r-r'} \lvert \ln \epsilon \rvert$, then we can distinguish two cases:
\begin{itemize}
\item (Non-resonant case)   if for any $n\in \Z^{d}$ with $0<\lvert n \rvert \leqslant N$, we have
$$
\lvert 2\rho - <n,\alpha> \rvert\geqslant \epsilon^{\sigma},
$$
then
$$\lvert B-Id\rvert_{r'}\leqslant \epsilon^{\frac{1}{2}} ,\   \ \lvert f_{+}\rvert_{r'}\leqslant 4\epsilon^{3-2\sigma}.$$
and
$$\lVert A_+-A\rVert\leqslant 2\lVert A\rVert\epsilon.$$
\item (Resonant case) if there exists $n_\ast$ with $0<\lvert n_\ast\rvert \leqslant N$ such that
$$
\lvert 2\rho- <n_\ast,\alpha> \rvert< \epsilon^{\sigma},
$$
then
$$\lvert B \rvert_{r'}\leqslant \epsilon^{-\frac{\sigma}{10}}\times\epsilon^{\frac{-r'}{r-r'}},\ \ \lVert B\rVert_0 \leqslant \epsilon^{-\frac{\sigma}{10}},\ \ \lvert f_{+}\rvert_{r'}\ll \epsilon^{100}.$$ Moreover, $A_+=e^{A''}$ with $\lVert A''\rVert \leqslant 2\epsilon^{\sigma}$.
\end{itemize}
\end{Proposition}

\begin{pf}
Although the result was stated for the $SL(2,\R)$ case, we prefer to prove it in the isomorphic group $SU(1,1)$, since the proof will become more explicit.  Notice that $sl(2,\R)$ is isomorphic to $su(1,1)$, which consists of matrices of the form
$$
\begin{pmatrix} it & v\\ \bar{v} & -it \end{pmatrix},
$$
with $t\in \R$, $v\in\C$.
The isomorphism between them is given by $A\rightarrow MAM^{-1}$, where
$$
M=\frac{1}{1+i}\begin{pmatrix} 1 & -i\\ 1 & i \end{pmatrix},
$$
and a simple calculation yields
$$
M\begin{pmatrix} x & y+z\\ y-z & -x \end{pmatrix}M^{-1}=\begin{pmatrix} iz & x-iy\\ x+iy & -iz \end{pmatrix}.
$$
Now we divide the proof into two cases:

\bigskip
 \textbf{Non-resonant case:}
For $0<\lvert n \rvert \leqslant N=\frac{2}{r-r'} \lvert \ln \epsilon \rvert$, we have
\begin{equation}\label{est1}
\lvert 2\rho - <n,\alpha> \rvert\geqslant \epsilon^{\sigma};
\end{equation}
by $(\ref{estf})$ with $D\geqslant 110$, we have
\begin{equation}\label{est2}
\left \lvert <n,\alpha> \right \rvert \geqslant \frac{\kappa}{\left \lvert n \right \rvert ^{\tau}}\geqslant \frac{\kappa}{\left \lvert N \right \rvert ^{\tau}}\geqslant \epsilon^{{\sigma}}.
\end{equation}
Define
\begin{equation}\label{lambdaN}
\Lambda_N=\{f\in C^{\omega}_{r}(\T^{d},su(1,1))\mid f(\theta)=\sum_{k\in \Z^{d},0<\lvert k \rvert<N}\hat{f}(k)e^{i<k,\theta>}\}.
\end{equation}
By $(\ref{est1})$ and $(\ref{est2})$, direct computation shows that if $Y\in\Lambda_N$, then
$$
\lvert A^{-1}Y(\theta+\alpha)A-Y(\theta)\rvert_r\geqslant\epsilon^{3\sigma}\lvert Y(\theta)\rvert_r,
$$
thus $\Lambda_N \subset \mathcal{B}_r^{nre}(\epsilon^{3\sigma})$.
Since $\epsilon^{3\sigma}\geqslant 13\lVert A\rVert^2\epsilon^{\frac{1}{2}}$, by Lemma \ref{lem2} we have $Y\in \mathcal{B}_r$ and $f^{re}\in \mathcal{B}_r^{re}(\epsilon^{3\sigma})$ such that
$$
e^{Y(\theta+\alpha)}(Ae^{f(\theta)})e^{-Y(\theta)}=Ae^{f^{re}(\theta)},
$$
with $\lvert Y \rvert_r\leqslant \epsilon^{\frac{1}{2}}$ and
\begin{equation}\label{fre}
\lvert f^{re}\rvert_r\leqslant 2\epsilon.
\end{equation}
By $(\ref{lambdaN})$
$$
(\mathcal{T}_N{f^{re}})(\theta)=\hat{f}^{re}(0), \ \ \lVert \hat{f}^{re}(0)\rVert \leqslant 2\epsilon,
$$
and
\begin{flalign}\label{13}
\lvert (\mathcal{R}_N{f^{re}})(\theta)\rvert_{r'}&= \lvert \sum_{\lvert n \rvert>N}\hat{f}^{re}(n)e^{i<n,\theta>}\rvert_{r'}\\
\notag &\leqslant 2\epsilon e^{-N(r-r')}(N)^{d}\\
\notag &\leqslant 2\epsilon^{3-2\sigma}.
\end{flalign}
Moreover, we can compute that
$$
e^{\hat{f}^{re}(0)+\mathcal{R}_N{f^{re}}(\theta)}=e^{\hat{f}^{re}(0)}(Id+e^{-\hat{f}^{re}(0)}\mathcal{O}(\mathcal{R}_N{f^{re}}))=e^{\hat{f}^{re}(0)}e^{f_+(\theta)},
$$
by $(\ref{13})$, we have
$$
\lvert f_+(\theta)\rvert_{r'}\leqslant 2\lvert \mathcal{R}_N{f^{re}(\theta)}\rvert_{r'} \leqslant 4\epsilon^{3-2\sigma}.
$$
Finally, if we denote
$$
A_+=Ae^{\hat{f}^{re}(0)},
$$
then we have
$$
\lVert A_+-A\rVert\leqslant \lVert A\rVert \lVert Id-e^{\hat{f}^{re}(0)} \rVert \leqslant 2\lVert A\rVert\epsilon.
$$

\bigskip
  \textbf{Resonant case:}
 Note that we only need to consider the case in which $A$ is elliptic with eigenvalues $\{e^{i\rho},e^{-i\rho}\}$ for $\rho\in \R\backslash\{0\}$, because if $\rho\in i\R$, the non-resonant condition is always fulfilled due to the Diophantine condition on $\alpha$ and then it indeed belongs to the non-resonant case.

\textbf{Claim:} $n_\ast$ is the unique resonant site with
$$
0<\lvert n_\ast\rvert \leqslant N=\frac{2}{r-r'} \lvert \ln \epsilon \rvert.
$$
\begin{pf}Indeed, if there exists  $n_{\ast}^{'}\neq n_\ast$ satisfying $|2\rho- <n_{\ast}^{'},\alpha>|< \epsilon^{\sigma}$, then by the Diophantine condition of $\alpha$, we have
$$
\frac{\kappa}{\lvert n_{\ast}^{'}-n_\ast\rvert^{\tau}}\leqslant \lvert <n_{\ast}^{'}-n_\ast,\alpha >\rvert< 2\epsilon^{\sigma},
$$
which implies that
$\lvert n_{\ast}^{'} \rvert>2^{-\frac{1}{\tau}}\kappa^{\frac{1}{\tau}}\epsilon^{-\frac{\sigma}{\tau}}-N\gg N.$\end{pf}

Since we have \begin{equation}\label{reso}
\lvert 2\rho- <n_\ast,\alpha> \rvert< \epsilon^{\sigma},
\end{equation}
 the smallness condition on $\epsilon$  implies that

$$\lvert \ln \epsilon\lvert ^\tau \epsilon^\sigma \leqslant \frac{\kappa (r-r')^\tau}{2^{\tau+1}}$$

\noindent Thus

$$\frac{\kappa}{\lvert n_\ast\lvert^\tau} \leqslant \lvert <n_\ast,\alpha>\lvert\leqslant \epsilon^\sigma +2\lvert \rho\lvert \leqslant \frac{\kappa}{2\lvert n_\ast\lvert ^\tau}+2\lvert \rho\lvert ,$$

\noindent which implies that

$$
\lvert \rho \rvert\geqslant \frac{\kappa}{4\lvert n_\ast\rvert^{\tau}}.
$$
Thus one can find $P\in SU(1,1)$ with
$$
\lVert P \rVert \leqslant \frac{2\lVert A \rVert}{\lvert \rho \rvert}\leqslant \frac{8\lVert A\rVert}{\kappa}\lvert n_\ast \rvert^{\tau},
$$
such that
$$
PAP^{-1}=\begin{pmatrix} e^{i\rho} & 0\\ 0 & e^{-i\rho} \end{pmatrix}=A'.
$$
Denote $g=PfP^{-1}$, by $(\ref{estf})$ with $D\geqslant 110$, we have
$$
\lVert A\rVert \lvert n_\ast\rvert^{\tau} \leqslant \lVert A\rVert N^{\tau}\leqslant \epsilon^{-\frac{\sigma}{11}},
$$
thus one gets the following estimates:
\begin{eqnarray}
\label{esti-p}\lVert P \rVert &\leqslant& \frac{8\lVert A\rVert}{\kappa}N^{\tau}\leqslant \frac{1}{2}\epsilon^{-\frac{\sigma}{10}},\\
\label{esti-g}\lvert g \rvert_r &\leqslant& \lVert P \rVert^2\lvert f\rvert_r \leqslant  \frac{1}{4}\epsilon^{1-\frac{\sigma}{5}}:=\epsilon'.
\end{eqnarray}
Now we define
\begin{flalign*}
&\Lambda_1(\epsilon^{\sigma})=\{n\in\Z^{d}: \lvert <n,\alpha>\rvert \geqslant \epsilon^{\sigma}\},\\
&\Lambda_2(\epsilon^{\sigma})=\{n\in\Z^{d}: \lvert 2\rho-<n,\alpha>\rvert \geqslant \epsilon^{\sigma}\},
\end{flalign*}
and define $\mathcal{B}_r^{nre}(\epsilon^{\sigma})$ as in $(\ref{space})$ with $A$ substituted by $A'$, then we can compute that any $Y\in \mathcal{B}_r^{nre}(\epsilon^{\sigma})$ takes the precise form:
\begin{small}
\begin{equation}
\begin{split}
Y(\theta)=&\sum_{n\in \Lambda_1(\epsilon^{\sigma})}\begin{pmatrix} i\hat{t}(n) & 0\\ 0 & -i\hat{t}(n) \end{pmatrix} e^{i<n,\theta>}+\\
&\sum_{n\in \Lambda_2(\epsilon^{\sigma})}\begin{pmatrix} 0 & \hat{v}(n)e^{i<n,\theta>}\\ \overline{\hat{v}(n)}e^{-i<n,\theta>} & 0 \end{pmatrix}.
\end{split}
\end{equation}
\end{small}
Since $\epsilon^{\sigma}\geqslant 13\lVert A'\rVert^2 (\epsilon') ^{\frac{1}{2}}$,  we can apply  Lemma \ref{lem2} to remove all the non-resonant terms of $g$, which means there exist $Y\in \mathcal{B}_r$ and $g^{re}\in \mathcal{B}_r^{re}(\eta)$ such that
$$
e^{Y(\theta+\alpha)}(A'e^{g(\theta)})e^{-Y(\theta)}=A'e^{g^{re}(\theta)},
$$
with $\lvert Y \rvert_r\leqslant (\epsilon')^{\frac{1}{2}}$ and $\lvert g^{re}\rvert_r\leqslant 2\epsilon'$.

Combining with the Diophantine condition on the frequency $\alpha$ and the Claim, we have:
\begin{flalign*}
&\{\Z^{d}\backslash\Lambda_1(\epsilon^{\sigma})\}\cap\{n\in \Z^{d}:\lvert n \rvert\leqslant \kappa^{\frac{1}{\tau}}\epsilon^{-\frac{\sigma}{\tau}}\}=\{0\},\\
&\{\Z^{d}\backslash\Lambda_2(\epsilon^{\sigma})\}\cap\{n\in \Z^{d}:\lvert n \rvert\leqslant 2^{-\frac{1}{\tau}}\kappa^{\frac{1}{\tau}}\epsilon^{-\frac{\sigma}{\tau}}-N\}=\{n_\ast\}.
\end{flalign*}
Let $N':=2^{-\frac{1}{\tau}}\kappa^{\frac{1}{\tau}}\epsilon^{-\frac{\sigma}{\tau}}-N$, then we can rewrite $g^{re}(\theta)$ as \begin{flalign*}
g^{re}(\theta)&=g^{re}_0+g^{re}_1(\theta)+g^{re}_2(\theta)\\
&=\begin{pmatrix} i\hat{t}(0) & 0 \\ 0 & -i\hat{t}(0) \end{pmatrix}+\begin{pmatrix} 0 & \hat{v}(n_\ast)e^{i<n_\ast,\theta>} \\ \overline{\hat{v}(n_\ast)}e^{-i<n_\ast,\theta>} & 0 \end{pmatrix}\\
&+\sum_{\lvert n \rvert>N'}\hat{g}^{re}(n)e^{i<n,\theta>}
\end{flalign*}


Define the $4\pi\Z$-periodic rotation $Q(\theta)$ as below:
$$
Q(\theta)=\begin{pmatrix} e^{-\frac{<n_\ast,\theta>}{2}i} & 0\\ 0 & e^{\frac{<n_\ast,\theta>}{2}i} \end{pmatrix}.
$$
So we have \begin{equation}\label{esti-Q}\lvert Q(\theta)\rvert_{r'}\leqslant e^{\frac{1}{2}Nr'}\leqslant \epsilon^{\frac{-r'}{r-r'}}.\end{equation} One can also show that
$$
Q(\theta+\alpha)(A'e^{g^{re}(\theta)})Q^{-1}(\theta)=\tilde{A}e^{\tilde{g}(\theta)},
$$
where
$$\tilde{A}=Q(\theta+\alpha)A'Q^{-1}(\theta)=\begin{pmatrix} e^{i(\rho-\frac{<n_\ast,\omega>}{2})} & 0\\ 0 & e^{-i(\rho-\frac{<n_\ast,\omega>}{2})} \end{pmatrix}$$
and
$$\tilde{g}(\theta)=Qg^{re}(\theta)Q^{-1}=Qg^{re}_0Q^{-1}+Qg^{re}_1(\theta)Q^{-1}+Qg^{re}_2(\theta)Q^{-1}.$$
Moreover,
\begin{flalign*}
Qg^{re}_0Q^{-1} &=g^{re}_0 = \begin{pmatrix} i\hat{t}(0) & 0 \\ 0 & -i\hat{t}(0) \end{pmatrix}
\in su(1,1), \\
Qg^{re}_1(\theta)Q^{-1}&=\begin{pmatrix} 0 & \hat{v}(n_\ast) \\ \overline{\hat{v}(n_\ast)} & 0 \end{pmatrix} \in su(1,1),
\end{flalign*}
Now we return back from $su(1,1)$ to $sl(2,\R)$. Denote
\begin{flalign*}
L &=M^{-1}(Qg^{re}_0Q^{-1}+Qg^{re}_1(\theta)Q^{-1})M,\\
F &=M^{-1}Qg^{re}_2(\theta)Q^{-1}M,\\
B &=M^{-1}(Q\circ e^Y \circ P) M,\\
\tilde{A}^{'} &=M^{-1}\tilde{A}M,
\end{flalign*}
then we have:
\begin{equation}\label{con1}
B(\theta+\alpha)(Ae^{f(\theta)})B^{-1}(\theta)=\tilde{A}^{'}e^{L+F(\theta)},
\end{equation}
by $(\ref{esti-p})$ and $(\ref{esti-Q})$, we have the following estimates:
\begin{eqnarray}
 \lVert B\rVert_0 &\leqslant&  |e^{Y}|_r \lVert P\rVert \leqslant \epsilon^{-\frac{\sigma}{10}},\\
\lvert B\rvert_{r'} &\leqslant& \epsilon^{-\frac{\sigma}{10}}\times\epsilon^{\frac{-r'}{r-r'}},\\
\label{D}\lVert L \rVert &\leqslant& \lVert Qg^{re}_0Q^{-1}\rVert + \lVert Qg^{re}_1(\theta)Q^{-1}\rVert \leqslant 2\epsilon^{1-\frac{\sigma}{5}},\\
\label{F}|F|_{r'} &\leqslant& \lvert Qg^{re}_2(\theta)Q^{-1}\rvert_{r'} \leqslant c\epsilon^{1-\frac{11\sigma}{5}}e^{-N'(r-r')}e^{Nr'}\leqslant \epsilon_{+}\ll \epsilon^{100}.
\end{eqnarray}

By $(\ref{D})$ and $(\ref{F})$, direct computation shows that
\begin{equation}\label{impl}
e^{L+F(\theta)}=e^L+\mathcal{O}(F(\theta))=e^L(Id+e^{-L}\mathcal{O}(F(\theta)))=e^L e^{f_+{(\theta)}}.
\end{equation}
It immediately implies that
$$
\lvert f_+{(\theta)}\rvert_{r'}\leqslant 2 |F(\theta)|_{r'}\leqslant 2 \epsilon_+ \ll \epsilon^{100}.
$$
Thus we can rewrite  $(\ref{con1})$ as
$$
B(\theta+\alpha)(Ae^{f(\theta)})B^{-1}(\theta)=A_{+}e^{f_+(\theta)},
$$
with
$$
A_+=\tilde{A}^{'}e^L=e^{A''}, \ \ A''\in sl(2,\R).
$$
Then $(\ref{reso})$ and $(\ref{D})$ gives
$$
\lVert A'' \rVert\leqslant 2(\lvert \rho-\frac{<n_\ast,\omega>}{2}\rvert+\lVert Qg^{re}_0Q^{-1}+Qg^{re}_1(\theta)Q^{-1}\rVert)\leqslant 2\epsilon^{\sigma}.
$$
This finishes the proof of Proposition $\ref{prop1}$.
\end{pf}

\subsection{Finitely differentiable case}

Now we are ready to transform the estimates in the analytic case into those in the $C^k$ case by analytic approximation.

\bigskip

Let $(f_{j})_{j\geqslant 1}$, $f_{j}\in C_{\frac{1}{j}}^{\omega}(\T^{d},sl(2,\R))$ be the analytic sequence approximating $f\in C^k(\T^{d},sl(2,\R))$, which satisfies $(\ref{2.1})$. Let $c,D$ be the constants defined in Proposition \ref{prop1}. Denote $$\epsilon_0^{'}(r,r')=\frac{c}{(2\lVert A\rVert)^D}(r-r')^{D\tau}.$$
Then for any $k\geqslant 5D\tau$, one can easily compute that

\begin{equation}\label{lowerboundk}\frac{c}{(2\lVert A\rVert)^D m^{\frac{k}{4}}}\leqslant \epsilon_0^{'}(\frac{1}{m},\frac{1}{m^2})
\end{equation}
for any $m\geqslant 10, m\in \Z$.

We start from $m_0=M> \max\{\frac{(2\lVert A\rVert)^D}{c},10\}$ and denote $l_j=M^{2^{j-1}}, j\geqslant 1$,
\begin{equation}\label{epsilonm}\epsilon_m=\frac{c}{(2\lVert A\rVert)^D m^{\frac{k}{4}}}.
\end{equation}

\begin{Proposition}\label{pro3}
Let $\alpha \in DC(\kappa,$ $\tau)$, $A\in SL(2,\R)$, $\sigma=\frac{1}{10}$, $f\in C^k(\T^{d}$,
$sl(2,\R))$ with $k\geqslant 5D\tau$ and $f_j$ be as above, there exists $\epsilon_0=\epsilon_0(\kappa,\tau,d,k,\lVert A\rVert)$ such that if $\lVert f\rVert_k\leqslant \epsilon_0$, then there exist $B_{l_j}\in C^{\omega}_{\frac{1}{l_{j+1}}}(2\T^{d},SL(2,\R))$, $A_{l_j}\in SL(2,\R)$, $f_{l_j}^{'}\in C^{\omega}_{\frac{1}{l_{j+1}}}(\T^{d},sl(2,\R))$, such that
$$
B_{l_j}(\theta+\alpha)(Ae^{f_{l_j}(\theta)})B^{-1}_{l_j}(\theta)=A_{l_j}e^{f_{l_j}^{'}(\theta)},
$$
with estimates

\begin{equation}\label{estimBlj}\lvert B_{l_j}(\theta)\rvert_{\frac{1}{l_{j+1}}}\leqslant \epsilon_{l_{j}}^{-\frac{2\sigma}{5}},\ \ \lVert B_{l_j}(\theta)\rVert_0\leqslant \epsilon_{l_{j}}^{-\frac{\sigma}{4}},
\end{equation}

$$
\lvert f_{l_j}^{'}(\theta)\rvert_{\frac{1}{l_{j+1}}}\leqslant \frac{1}{2}\epsilon_{l_{j}}^{\frac{5}{2}},\ \ \lVert A_{l_j}\rVert\leqslant 2\lVert A\rVert.
$$
Moreover, there exists unitary $U_j\in SL(2,\C)$ such that
$$
U_jA_{l_j}U_j^{-1}=\begin{pmatrix} e^{\gamma_j} & c_j\\ 0 & e^{-\gamma_j} \end{pmatrix},
$$
and
\begin{equation}
\label{estisharp}\lvert c_j\rvert \lVert B_{l_j}(\theta)\rVert_0^2\leqslant 8\lVert A\rVert,
\end{equation}
with $\gamma_j\in i\R\cup\R$ and $c_j\in \C$.
\end{Proposition}

\begin{pf} \textbf{First step:} Assume that $$C'\lVert f(\theta) \rVert_k\leqslant \frac{c}{(2\lVert A\rVert)^D l_{1}^{\frac{k}{4}}},$$ by $(\ref{2.1})$ we have
$$
\lvert f_{l_{1}}(\theta)\rvert_{\frac{1}{l_{1}}}\leqslant \epsilon_{l_{1}}\leqslant \epsilon_0^{'}(\frac{1}{l_{1}},\frac{1}{l_{2}}).
$$
Apply Proposition \ref{prop1}, we can find $B_{l_1}\in C^{\omega}_{\frac{1}{l_2}}(2\T^{d},SL(2,\R))$, $A_{l_1}\in SL(2,\R)$ and $f_{l_1}^{'}\in C^{\omega}_{\frac{1}{l_2}}(\T^{d},sl(2,\R))$ such that
$$
B_{l_1}(\theta+\alpha)(Ae^{f_{l_1}(\theta)})B^{-1}_{l_1}(\theta)=A_{l_1}e^{f_{l_1}^{'}(\theta)},
$$
More precisely, we have two different cases:
\begin{itemize}
\item (Non-resonant case)
$$
\lvert B_{l_1}\rvert_{\frac{1}{l_2}}\leqslant 1+\epsilon_{l_{1}}^{\frac{1}{2}}, \ \ \lvert f_{l_1}^{'}\rvert_{\frac{1}{l_2}}\leqslant 4\epsilon_{l_{1}}^{3-2\sigma}\leqslant  \frac{1}{2}\epsilon_{l_{1}}^{\frac{5}{2}},
$$
and
$$
\lVert A_{l_1}-A\rVert \leqslant 2\lVert A\rVert\epsilon_{l_{1}}.
$$
\item (Resonant case)
$$
\lvert B_{l_1}\rvert_{\frac{1}{l_2}}\leqslant \epsilon_{l_{1}}^{-\frac{\sigma}{10}}\times\epsilon_{l_{1}}^{\frac{-1}{l_1-1}}\leqslant \epsilon_{l_{1}}^{-\frac{\sigma}{5}},
\ \ \lvert B_{l_1}\rvert_0 \leqslant \epsilon_{l_{1}}^{-\frac{\sigma}{10}},
$$
$$
\lvert f_{l_1}^{'}\rvert_{\frac{1}{l_2}}\ll \epsilon_{l_{1}}^{100} \leqslant \frac{1}{2}\epsilon_{l_{1}}^{\frac{5}{2}}.
$$
Moreover, $A_{l_1}=e^{A_{l_1}''}$ with $\lVert A_{l_1}''\rVert \leqslant 2\epsilon_{l_{1}}^{\sigma}$.
\end{itemize}

\textbf{Induction step:} Assume that for $l_n,n\leqslant k$, we already have $(\ref{estisharp})$ and
\begin{equation}
\label{estiln}B_{l_n}(\theta+\alpha)(Ae^{f_{l_n}(\theta)})B^{-1}_{l_n}(\theta)=A_{l_n}e^{f_{l_n}^{'}(\theta)},
\end{equation}
with
\begin{equation}
\label{estima}\lvert B_{l_n}(\theta)\rvert_{\frac{1}{l_{n+1}}}\leqslant \epsilon_{l_n}^{-\frac{2\sigma}{5}},\ \ \lVert B_{l_n}(\theta)\rVert_0\leqslant \epsilon_{l_n}^{-\frac{\sigma}{4}},\ \ \lvert f_{l_n}^{'}(\theta)\rvert_{\frac{1}{l_{n+1}}}\leqslant \frac{1}{2}\epsilon_{l_n}^{\frac{5}{2}},
\end{equation}
and
\begin{equation}
\label{esta}\lVert A_{l_n}\rVert \leqslant 2\lVert A\rVert.
\end{equation}
Moreover, if the $n$-th step is obtained by the resonant case, we have
\begin{equation}
\label{estres}A_{l_n}=e^{A_{l_n}''}, \ \ \lVert A_{l_n}''\rVert <2\epsilon_{l_{n}}^{\sigma}.
\end{equation}
If the $n$-th step is obtained by the non-resonant case, we have
\begin{equation}
\label{estnon}\lVert A_{l_n}-A_{l_{n-1}}\rVert\leqslant 2\lVert A_{l_{n-1}}\rVert \epsilon_{l_{n}}.
\end{equation}
and
\begin{equation}
\label{estb}\lVert B_{l_n}\rVert_{\frac{1}{l_{n+1}}}\leqslant (1+\epsilon_{l_{n}}^{\frac{1}{2}})\lVert B_{l_{n-1}}\rVert_{\frac{1}{l_{n}}}.
\end{equation}
Now by $(\ref{estiln})$, for $l_{n+1}, n=k$, we have
$$
B_{l_n}(\theta+\alpha)(Ae^{f_{l_{n+1}}})B^{-1}_{l_n}(\theta)=A_{l_n}e^{f_{l_n}^{'}}+B_{l_n}(\theta+\alpha)(Ae^{f_{l_{n+1}}}-Ae^{f_{l_{n}}})B^{-1}_{l_n}(\theta).
$$
In $(\ref{2.1})$, a simple integral implies
\begin{equation}\label{ov}
\lvert f_{l_{n+1}}(\theta)-f_{l_{n}}(\theta)\rvert_{\frac{1}{l_{n+1}}}\leqslant \frac{c}{(2\lVert A\rVert)^D l_{n}^{k-1}},
\end{equation}
Moreover, $(\ref{2.1})$ also gives us
\begin{equation}\label{ok}
\lvert f_{l_{n+1}}(\theta)\rvert_{\frac{1}{l_{n+1}}} +\lvert f_{l_{n}}(\theta)\rvert_{\frac{1}{l_{n+1}}}\leqslant \frac{2c}{(2\lVert A\rVert)^D M^{\frac{k}{4}}}.
\end{equation}
Thus if we rewrite that $$A_{l_n}e^{f_{l_n}^{'}(\theta)}+B_{l_n}(\theta+\alpha)(Ae^{f_{l_{n+1}}(\theta)}-Ae^{f_{l_{n}}(\theta)})B^{-1}_{l_n}(\theta)=A_{l_n}e^{\widetilde{f_{l_n}}(\theta)},$$
by $(\ref{estima})$, $(\ref{esta})$, $(\ref{ov})$ and $(\ref{ok})$ we obtain
\begin{flalign*}
\lvert \widetilde{f_{l_n}}(\theta)\rvert_{\frac{1}{l_{n+1}}}\leqslant & \lvert f_{l_n}^{'}(\theta)\rvert_{\frac{1}{l_{n+1}}}+\lVert A_{l_n}^{-1}\rVert\lvert B_{l_n}(\theta+\alpha)(Ae^{f_{l_{n+1}}(\theta)}-Ae^{f_{l_{n}}(\theta)})B^{-1}_{l_n}(\theta)\rvert_{\frac{1}{l_{n+1}}}\\
\leqslant &\frac{1}{2}\epsilon_{l_{n}}^{\frac{5}{2}}+ 2\lVert A\rVert^2\times\epsilon_{l_n}^{-\frac{4\sigma}{5}}\times\frac{c}{(2\lVert A\rVert)^Dl_{n}^{k-1}}\\
\leqslant & \frac{1}{2}\epsilon_{l_{n+1}}+\frac{1}{2}\times\frac{c}{(2\lVert A\rVert)^Dl_{n}^{\frac{k}{2}}}\\
\leqslant & \epsilon_{l_{n+1}}\\
\leqslant & \epsilon_0^{'}(\frac{1}{l_{n+1}},\frac{1}{l_{n+2}}).
\end{flalign*}
Now for $(\alpha,A_{l_n}e^{\widetilde{f_{l_n}}(\theta)})$, we can apply Proposition \ref{prop1} again to get $\tilde{B}_{l_n}\in C^{\omega}_{\frac{1}{l_{n+2}}}(2\T^{d},$
$SL(2,\R))$, $A_{l_{n+1}}\in SL(2,\R)$ and $f_{l_{n+1}}^{'}\in C^{\omega}_{\frac{1}{l_{n+2}}}(\T^{d},sl(2,\R))$ such that
$$
\tilde{B}_{l_n}(\theta+\alpha)(A_{l_n}e^{\widetilde{f_{l_n}}(\theta)})\tilde{B}_{l_n}(\theta)=A_{l_{n+1}}e^{f_{l_{n+1}}^{'}(\theta)},
$$
with
$$\lvert\tilde{B}_{l_n}(\theta)\rvert_{\frac{1}{l_{n+2}}}\leqslant \epsilon_{l_{n+1}}^{-\frac{\sigma}{10}}\times\epsilon_{l_{n+1}}^{\frac{-1}{l_{n+1}-1}},\ \ \lVert \tilde{B}_{l_n}(\theta)\rVert_0\leqslant \epsilon_{l_{n+1}}^{-\frac{\sigma}{10}},\ \ \lvert f_{l_{n+1}}^{'}(\theta)\rvert_{\frac{1}{l_{n+2}}}\leqslant \frac{1}{2}\epsilon_{l_{n+1}}^{\frac{5}{2}}.$$

Denote $B_{l_{n+1}}:=\tilde{B}_{l_n}B_{l_{n}} \in C^{\omega}_{\frac{1}{l_{n+2}}}(2\T^{d},SL(2,\R))$, then
$$
\lvert B_{l_{n+1}}(\theta)\rvert_{\frac{1}{l_{n+2}}}\leqslant \epsilon_{l_n}^{-\frac{2\sigma}{5}}\times\epsilon_{l_{n+1}}^{-\frac{\sigma}{10}}\times\epsilon_{l_{n+1}}^{\frac{-1}{l_{n+1}-1}}\leqslant \epsilon_{l_{n+1}}^{-\frac{2\sigma}{5}},
$$
and
$$
\lVert B_{l_{n+1}}(\theta) \rVert_0\leqslant \epsilon_{l_n}^{-\frac{\sigma}{4}}\times \epsilon_{l_{n+1}}^{-\frac{\sigma}{10}} \leqslant \epsilon_{l_{n+1}}^{-\frac{\sigma}{4}}.
$$
For the remaining estimates, we distinguish two cases.

\bigskip
$\bullet$ If the $(n+1)$-th step is in the resonant case, we have
$$
A_{l_{n+1}}=e^{A_{l_{n+1}}''}, \ \ \lVert A_{l_{n+1}}''\rVert < 2\epsilon_{l_{{n+1}}}^{\sigma}, \ \ \lVert A_{l_{n+1}}\rVert\leqslant 1+2\epsilon_{l_{{n+1}}}^{\sigma}\leqslant 2\lVert A\rVert.
$$
Then there exists unitary $U\in SL(2,\C)$ such that
\begin{equation}
\label{estt}UA_{l_{n+1}}U^{-1}=\begin{pmatrix} e^{\gamma_{n+1}} & c_{n+1}\\ 0 & e^{-\gamma_{n+1}} \end{pmatrix},
\end{equation}
with $\lvert c_{n+1}\rvert\leqslant 2\lVert A_{l_{n+1}}''\rVert \leqslant 4\epsilon_{l_{{n+1}}}^{\sigma}$. Thus $(\ref{estisharp})$ is fulfilled.

\bigskip
$\bullet$
If the $(n+1)$-th step is in the non-resonant case, one traces back to the resonant step $j$ which is closest to $n+1$.

If $j$ exists, by $(\ref{estima})$ and $(\ref{estres})$ we have
$$
\lvert B_{l_j}(\theta)\rvert_{\frac{1}{l_{j+1}}}\leqslant \epsilon_{l_j}^{-\frac{2\sigma}{5}},\ \ \lVert B_{l_j}(\theta)\rVert_0\leqslant \epsilon_{l_j}^{-\frac{\sigma}{4}},
$$
$$
A_{l_j}=e^{A_{l_j}''}, \ \ \lVert A_{l_j}''\rVert < 2\epsilon_{l_{j}}^{\sigma}, \ \ \lVert A_{l_j}\rVert \leqslant 1+2\epsilon_{l_{j}}^{\sigma}.
$$
By our choice of $j$, from $j$ to $n+1$, every step is non-resonant. Thus by $(\ref{estnon})$ we obtain
\begin{equation}
\label{estnn}\lVert A_{l_{n+1}}- A_{l_j}\rVert\leqslant 2\epsilon_{l_j}^{\frac{1}{2}},
\end{equation}
so
$$
\lVert A_{l_{n+1}}\rVert \leqslant 1+2\epsilon_{l_{j}}^{\sigma}+ 2\epsilon_{l_j}^{\frac{1}{2}}\leqslant 2\lVert A\rVert.
$$
Estimate $(\ref{estnn})$ implies that if we rewrite $A_{l_{n+1}}=e^{A_{l_{n+1}}''}$, then
$$
\lVert A_{l_{n+1}}''\rVert \leqslant 4\epsilon_{l_{j}}^{\sigma}.
$$
Moreover, by $(\ref{estb})$, we have
$$
\lVert B_{l_{n+1}}(\theta)\rVert_0 \leqslant \lVert B_{l_{n+1}}(\theta)\rVert_{\frac{1}{l_{n+1}}}\leqslant 2\lvert B_{l_j}(\theta)\rvert_{\frac{1}{l_{j+1}}}\leqslant 2\epsilon_{l_j}^{-\frac{2\sigma}{5}}.
$$
Similarly to the process of $(\ref{estt})$, $(\ref{estisharp})$ is fulfilled.

If $j$ vanishes, it immediately implies that from $1$ to $n+1$, each step is non-resonant. In this case, $\lVert A_{l_{n+1}}\rVert \leqslant 2\lVert A\rVert$ and the estimate $(\ref{estisharp})$ is naturally satisfied as
$$\lVert B_{l_{n+1}}(\theta)\rVert_{\frac{1}{l_{n+1}}}\leqslant 2.$$
\end{pf}

\subsection{Differentiable almost reducibility}

Now, we are ready to show the quantitative almost reducibility for $C^{k}$ quasi-periodic $SL(2,\R)$ cocycles.

\begin{Theorem}\label{thm3}
Let $\alpha\in DC(\kappa,\tau)$, $A\in SL(2,\R)$, $\sigma=\frac{1}{10}$, $f\in C^k(\T^{d},sl(2,\R))$ with $k\geqslant 5D\tau$, there exists $\epsilon_1=\epsilon_1(\kappa,\tau,d,k,\lVert A\rVert)$ such that if $\lVert f\rVert_k\leqslant \epsilon_1$ then $(\alpha , Ae^{f(\theta)})$ is $C^{k,k_0}$ almost reducible with $k_0\in\N$, $k_0\leqslant \frac{k}{6}$. Moreover, if we further assume
$(\alpha , Ae^{f(\theta)})$ is not uniformly hyperbolic, then there exists $B_{l_j}\in C^{\omega}_{\frac{1}{l_{j+1}}}(2\T^{d}$,
$SL(2,\C))$, $A_{l_j}\in SL(2,\C)$, $\tilde{F}_{l_j}\in C^{k}(\T^{d}$,
$SL(2,\C))$, such that
$$
B_{l_j}(\theta+\alpha)(Ae^{f(\theta)})B^{-1}_{l_j}(\theta)=A_{l_j}+\tilde{F}_{l_j}(\theta)
$$
with $$\lVert B_{l_j}\rVert_0\leqslant \epsilon_{l_j}^{-\frac{\sigma}{4}},\ \ \lVert \tilde{F}_{l_j}\rVert_0\leqslant \epsilon_{l_j}^{\frac{1}{4}}$$
and
$A_{l_j}=\begin{pmatrix} e^{i\gamma_j} & c_j\\ 0 & e^{-i\gamma_j} \end{pmatrix}$ with estimate

\begin{equation}\label{ess}
\lVert B_{l_j}\rVert_0^2 \lvert c_j\rvert\leqslant 8\lVert A\rVert,
\end{equation}
 where $\gamma_j\in \R$ and $c_j \in \C$.
\end{Theorem}

\begin{Remark}\label{remark3} Estimate $(\ref{ess})$ is essential for us to obtain $1/2$-H\"older continuity.
\end{Remark}
\begin{pf}
We first deal with the $C^0$ estimate as it is much more significant for the proof of our main Theorem. By Proposition $\ref{pro3}$, we have for any $l_j$, $j\in \N^+$:
$$
B_{l_j}(\theta+\alpha)(Ae^{f_{l_j}(\theta)})B^{-1}_{l_j}(\theta)=A_{l_j}e^{f_{l_j}^{'}(\theta)},
$$
thus
$$
B_{l_j}(\theta+\alpha)(Ae^{f(\theta)})B^{-1}_{l_j}(\theta)=A_{l_j}e^{f_{l_j}^{'}(\theta)}+B_{l_j}(\theta+\alpha)(Ae^{f(\theta)}-Ae^{f_{l_j}(\theta)})B^{-1}_{l_j}(\theta).
$$
Denote
\begin{equation}\label{quat1}
A_{l_j}+\tilde{F}_{l_j}(\theta)=A_{l_j}e^{f_{l_j}^{'}(\theta)}+B_{l_j}(\theta+\alpha)(Ae^{f(\theta)}-Ae^{f_{l_j}(\theta)})B^{-1}_{l_j}(\theta).
\end{equation}
In $(\ref{2.1})$, by a simple integration we get
\begin{equation}\label{quat2}
\lVert f(\theta)-f_{l_j}(\theta)\rVert_0\leqslant \frac{c}{(2\lVert A\rVert)^Dl_j^{k-1}},
\end{equation}
and
\begin{equation}\label{quat3}
\lVert f(\theta)\rVert_0+\lVert f_{l_j}(\theta)\rVert_0\leqslant \frac{c}{(2\lVert A\rVert)^DM^{\frac{k}{4}}}+\frac{c}{(2\lVert A\rVert)^DC'M^{\frac{k}{4}}}.
\end{equation}
Proposition $\ref{pro3}$ also gives the estimates
\begin{equation}\label{quat4}
\lVert B_{l_j}(\theta)\rVert_0\leqslant \epsilon_{l_j}^{-\frac{\sigma}{4}},\ \ \lvert f_{l_j}^{'}(\theta)\rvert_{\frac{1}{l_{j+1}}}\leqslant\frac{1}{2}\epsilon_{l_j}^{\frac{5}{2}},
\end{equation}
and
\begin{equation}\label{quat5}
\lVert A_{l_j}\rVert\leqslant 2\lVert A\rVert.
\end{equation}
Thus by $(\ref{quat1})-(\ref{quat5})$, we have
\begin{flalign*}
\lVert \tilde{F}_{l_j}(\theta)\rVert_0 \leqslant & \lVert A_{l_j}f_{l_j}^{'}(\theta)\rVert_0+\lVert B_{l_j}(\theta+\alpha)(Ae^{f(\theta)}-Ae^{f_{l_j}(\theta)})B^{-1}_{l_j}(\theta)\rVert_0\\
\leqslant &\lVert A\rVert\epsilon_{l_j}^{\frac{5}{2}}+\lVert A\rVert\times\epsilon_{l_j}^{-\frac{\sigma}{2}}\times\frac{c}{(2\lVert A\rVert)^D l_j^{k-1}} \leqslant \frac{\epsilon_{l_j}^2}{2\lVert A \rVert}.
\end{flalign*}

To prove $(\ref{ess})$, Proposition \ref{pro3} implies that we only need to exclude the case when $\gamma_j\in i\R \backslash \{0\}$.

We assume that $spec(A_{l_j})=\{e^{\lambda_j},e^{-\lambda_j}\}, \lambda_j\in \R \backslash \{0\}$, thus we can always find $P\in SO(2,\R)$ such that
$$
PA_{l_j}P^{-1}=\begin{pmatrix} e^{\lambda_j} & c_j \\ 0 & e^{-\lambda_j}\end{pmatrix},
$$
with $\lvert c_j \rvert \leqslant \lVert A_{l_j}\rVert\leqslant 2\lVert A\rVert$.

If $\lvert\lambda_j\rvert>\epsilon_{l_j}^{\frac{1}{4}}$, set $B=diag\{ \lVert 2A\rVert^{-\frac{1}{2}}\epsilon_{l_j}^{\frac{1}{2}}, \lVert 2A\rVert^{\frac{1}{2}} \epsilon_{l_j}^{-\frac{1}{2}} \}$, then
\begin{equation}\label{123}
BP(A_{l_j}+\tilde{F}_{l_j}(\theta))P^{-1}B^{-1}=\begin{pmatrix} e^{\lambda_j} & 0 \\ 0 & e^{-\lambda_j}\end{pmatrix}+ F(\theta),
\end{equation}
where $\lVert F(\theta)\rVert_0 \leqslant 2\epsilon_{l_j}$.
We rewrite
$$
\begin{pmatrix} e^{\lambda_j} & 0 \\ 0 & e^{-\lambda_j}\end{pmatrix}+ F(\theta)=\begin{pmatrix} e^{\lambda_j} & 0 \\ 0 & e^{-\lambda_j}\end{pmatrix}e^{\tilde{f}(\theta)}
$$
with $\lVert \tilde{f}(\theta)\rVert_0\leqslant 4\epsilon_{l_j}$.
Then by Remark $\ref{rem2}$ and Corollary 3.1 of \cite{houyou}, one can conjugate $(\ref{123})$ to
$$
\begin{pmatrix} e^{\lambda_j} & 0 \\ 0 & e^{-\lambda_j}\end{pmatrix}\begin{pmatrix} e^{\tilde{f}^{re}(\theta)} & 0 \\ 0 & e^{-\tilde{f}^{re}(\theta)}\end{pmatrix}=\begin{pmatrix} e^{\lambda_j}e^{\tilde{f}^{re}(\theta)} & 0 \\ 0 & e^{-\lambda_j}e^{-\tilde{f}^{re}(\theta)}\end{pmatrix}
$$
with $\lVert \tilde{f}^{re}(\theta)\rVert_0 \leqslant 8\epsilon_{l_j}$,
thus $(\alpha , Ae^{f(\theta)})$ is uniformly hyperbolic, which contradicts our assumption. Thus we only need to consider $\lvert\lambda_j\rvert \leqslant \epsilon_{l_j}^{\frac{1}{4}}$; in this case, we put $\lambda_j$ into the perturbation so that the new perturbation satisfies $\lVert \tilde{F}^{'}_{l_j} \rVert_0\leqslant \epsilon_{l_j}^{\frac{1}{4}}$ and
$$
A_{l_j}=\begin{pmatrix} 1 & c_j \\ 0 & 1\end{pmatrix}.
$$

Although $C^0$ norm is sufficient to prove our main Theorem, we want to show the actually differentiable almost reducibility instead of $C^0$ almost reducibility. By Cauchy estimates, for $k_0 \in \N$ and $k_0\leqslant k $, we have
\begin{flalign*}
\lVert f_j-f_{j+1}\rVert_{k_0}& \leqslant \sup_{\substack{
                             l\leqslant k_0,
                             \theta \in \T^{d}
                          }}\lVert (\partial_{\theta_1}^{l_1}+\cdots+\partial_{\theta_d}^{l_d})(f_j(\theta)-f_{j+1}(\theta))\rVert \\
                          &\leqslant (k_0)!(j+1)^{k_0}\lvert f_j-f_{j+1}\rvert_{\frac{1}{j+1}}\\
                          &\leqslant (k_0)!(j+1)^{k_0}\times \frac{c}{(2\lVert A\rVert)^D j^{k}}\\
                          &\leqslant \frac{C_1}{j^{k-k_0}}
\end{flalign*}
where $C_1$ is independent of $j$.

By a simple integration we get
$$
\lVert f(\theta)-f_{l_j}(\theta)\rVert_{k_0}\leqslant \frac{C_1}{l_j^{k-k_0-1}}.
$$
Similarly by Cauchy estimates, we have
\begin{flalign*}
\lVert f_{l_j}^{'}(\theta)\rVert_{k_0}& \leqslant (k_0)!(l_{j+1})^{k_0}\lvert f_{l_j}^{'}(\theta)\rvert_{\frac{1}{l_{j+1}}}\\
&\leqslant (k_0)!(l_{j})^{2k_0}\times \frac{1}{2}\epsilon_{l_j}^{\frac{5}{2}}\\
&\leqslant (k_0)!(l_{j})^{2k_0}\times \frac{1}{2}\times (\frac{c}{(2\lVert A\rVert)^D {l_j}^{\frac{k}{4}}})^{\frac{5}{2}}\\
&\leqslant \frac{C_2}{l_j^{\frac{5k}{8}-2k_0}}
\end{flalign*}
where $C_2$ is independent of $j$.
\begin{flalign*}
\lVert B_{l_j}(\theta)\rVert_{k_0}&\leqslant (k_0)!(l_{j+1})^{k_0}\lvert B_{l_j}(\theta)\rvert_{\frac{1}{l_{j+1}}}\\
&\leqslant (k_0)!(l_{j})^{2k_0}\times \epsilon_{l_{j}}^{-\frac{2\sigma}{5}}\\
&\leqslant (k_0)!(l_{j})^{2k_0}\times (\frac{c}{(2\lVert A\rVert)^D {l_j}^{\frac{k}{4}}})^{-\frac{2\sigma}{5}}\\
&\leqslant C_3 l_j^{\frac{k}{100}+2k_0}
\end{flalign*}
where $C_3$ is independent of $j$.
Thus we have
\begin{flalign*}
\lVert \tilde{F}_{l_j}(\theta)\rVert_{k_0} &\leqslant  \lVert A_{l_j}f_{l_j}^{'}(\theta)\rVert_{k_0}+\lVert B_{l_j}(\theta+\alpha)(Ae^{f(\theta)}-Ae^{f_{l_j}(\theta)})B^{-1}_{l_j}(\theta)\rVert_{k_0}\\
&\leqslant \frac{C_4}{l_j^{\frac{5k}{8}-2k_0}}+C_5 l_j^{\frac{k}{50}+4k_0}\times l_j^{-k+k_0+1}\\
&\leqslant \frac{C_4}{l_j^{\frac{5k}{8}-2k_0}}+\frac{C_5}{l_j^{\frac{49k}{50}-5k_0-1}}.
\end{flalign*}
If $k_0\leqslant \frac{k}{6}$, then
$$
\lVert \tilde{F}_{l_j}(\theta)\rVert_{k_0}\leqslant \frac{C_6}{l_j^2},
$$
which immediately shows that $\{ \tilde{F}_{l_j}\}_{j\geqslant 1}$ is a Cauchy sequence in $C^{k_0}$ topology with
$$\lim_{j \to +\infty}\lVert \tilde{F}_{l_j}(\theta)\rVert_{k_0}=0,$$
which means precisely that $(\alpha , Ae^{f(\theta)})$ is $C^{k,k_0}$ almost reducible.

This finishes the proof of Theorem \ref{thm3}.
\end{pf}

\subsection{Remarks on differentiable almost reducibility}\label{3.5}

The statement on differentiable almost reducibility was given in \cite{chavaudret3}, based on the almost reducibility result in \cite{chavaudret2}. However, the proof in \cite{chavaudret3} is not precise since it was based on  estimates on $B_j$ given in \cite{chavaudret2} which are not good enough. In fact, from the estimates in \cite{chavaudret2}, one gets the following estimates for the change of variables $B_j$ (as is usual in the analytic case):

$$\lvert B_j\lvert_{r_j}\leq (\epsilon_0^{-2^j})^{\frac{1}{4}(r_{j-1}-r_j)}
$$
\noindent where $r_j$ are the successive radii of analyticity. However, since  the perturbation comes from the analytic approximation $(\ref{2.1})$, which is of size $\frac{1}{r_{j-1}^k}$, when $r_j$ decreases as $\frac{1}{j}$ as in \cite{chavaudret2}. Thus  the bounds for the new perturbation is $$(\epsilon_0^{-2^j})^{\frac{1}{4}(r_{j-1}-r_j)} \frac{1}{r_{j-1}^k},$$ which is divergent, and thus the iteration can not go on.  For this reason, in this paper, we have to work out for $B_j$ a  better  estimates
$$(\epsilon_{j-1})^{\frac{-r_j}{r_{j-1}-r_j}},$$
which is enough for our purposes.  

Another improvement is the lower bound on the degree of differentiability $k$ (we obtain an integer $k$ which only depends on $\tau$, which fits the intuition from classical KAM theory), motivated by the estimate \eqref{lowerboundk} above, which is more explicit than the one in \cite{chavaudret3}.

However, we do not claim optimality of the condition we impose on $k$ and $k_0$. In fact, our KAM scheme still has quite a lot of flexibility.  As we pointed out, the new perturbation is of size $\frac{1}{r_{j-1}^k}$, where we have set the sequence $r_j=r_{j-1}^2$, which is more suitable for the resonant case. This has the advantage of providing a unified argument, but is clearly not optimal in the non-resonant case, since in this case, the conjugacy map is close to identity, which allows slow shrinking of the analytic radius (for example, from $1/j$ to $1/(j+1)$. More precise estimates should be obtainable by interpolating these two different arguments. These further modifications  will have two advantages: the first is that one can obtain $C^{k,k-D}$ almost reducibility results with $D$ independent of $k$; the second advantage is that one can obtain quite better estimates on the lower bound on the degree of differentiability $k$. However, since the main point of this paper is to obtain the H\"older continuity of LE, this rough estimate is already enough for our purposes.

\section{Proof of the main results}

As we mentioned in the introduction, nice quantitative almost reducibility of $(\alpha, A)$ easily implies H\"older continuity of the Lyapunov exponent, where \textquotedblleft nice\textquotedblright means that
the norm of the conjugation map is well controlled compared with the upper triangular element of the constant matrix and the norm of the perturbation. In fact, for $B$ close to $A$, we can write that $(\alpha, B)=(\alpha, A+(B-A))$, so nice quantitative almost reducibility of $(\alpha, A)$ implies nice quantitative almost reducibility of $(\alpha, B)$, which indicates the closeness between $L(\alpha, B)$ and $L(\alpha, A)$, as we shall see below.
%
%
\subsection{Proof of Theorem 1.1}
\begin{pf}
If the cocycle is uniformly hyperbolic, then the Lyapunov exponent is Lipschitz, which follows from the cone-field criterion. Therefore, we only need to consider the case $(\alpha,\mathcal{A})$ is almost reducible, but not uniformly hyperbolic.

 By the definition of $C^{k^{'},k}$ almost reducibility, there exist $B_j\in C^{k}(\T^{d},$ $SL(2,\R))$, $A\in SL(2,\R)$ and $F_j\in C^{k}(\T^{d},sl(2,\R))$ such that $$B_j(\theta+\alpha)\mathcal{A}(\theta)B_j(\theta)^{-1}=Ae^{F_j(\theta)}$$ with $\lVert F_j \rVert_{k}\leqslant \epsilon_j \rightarrow 0$. Furthermore, without loss of generality, we may assume $\|A\|\leqslant 1$.
%
%
%
%

Now we have $\lVert A\rVert \leqslant 1$ and $k\geqslant D_0\tau = 5D\tau$.
In the following, $C_0$ is a large constant depending on $d,\kappa,\tau$ and $c$ is a small constant depending on $d,\kappa,\tau$.
Notice that $(\alpha,Ae^{F_j(\theta)})$ is not uniformly hyperbolic, so for $j$ large enough, it will finally fall into our local regime in Theorem $\ref{thm3}$. By Theorem $\ref{thm3}$, we rewrite
$$
A_{l_j}+\tilde{F}_{l_j}(\theta)=\begin{pmatrix} e^{i\gamma_j} & 0 \\ 0 & e^{-i\gamma_j} \end{pmatrix}+\begin{pmatrix} q_1(\theta) & q_2(\theta) \\ q_3(\theta) & q_4(\theta) \end{pmatrix},
$$
where
$$
\lVert q_2(\theta)\rVert _0 \leqslant \lvert c_j\rvert+ \lVert \tilde{F}_{l_j}\rVert_0\leqslant \lvert c_j\rvert+\epsilon_{l_j}^{\frac{1}{4}},
$$
$$
\lVert q_1(\theta)\rVert _0,\lVert q_3(\theta)\rVert _0,\lVert q_4(\theta)\rVert _0 \leqslant 2\epsilon_{l_j}^{\frac{1}{4}}\leqslant \frac{C_0}{l_j^{\frac{k}{16}}}.
$$
Let $$D=\begin{pmatrix} d & 0 \\ 0 & d^{-1} \end{pmatrix}$$ where $d=\lVert B_{l_j}\rVert_0\epsilon^{\frac{1}{4}}$. Let $W(\theta)=DB_{l_j}(\theta)$. If $d \leqslant 1$, i.e. $$\epsilon \leqslant (\epsilon_{l_j}^{-\frac{\sigma}{4}})^{-4},$$ thus we pick $$\epsilon\leqslant c{l_j^{-\frac{k}{40}}},$$ then we have $$\lVert W \rVert_0\leqslant C_0\epsilon^{-\frac{1}{4}}.$$ Since
$$
D\begin{pmatrix} x_1 & x_2 \\ x_3 & x_4 \end{pmatrix}D^{-1}=\begin{pmatrix} x_1 & d^2 x_2 \\ d^{-2}x_3 & x_4 \end{pmatrix},
$$
we have
$$
Z(\theta)=W(\theta+\alpha)(Ae^f)W(\theta)^{-1}=\begin{pmatrix} e^{i\gamma_j} & 0 \\ 0 & e^{-i\gamma_j} \end{pmatrix}+\begin{pmatrix} q_1(\theta) & d^2 q_2(\theta)\\  d^{-2}q_3(\theta) & q_4(\theta) \end{pmatrix},
$$
with
$$\lVert q_1\rVert_0,\lVert q_4\rVert_0,\lVert d^{-2}q_3\rVert_0\leqslant C_0\epsilon^{-\frac{1}{2}}l_j^{-\frac{k}{16}},$$
and by $(\ref{ess})$ we have
$$\lVert d^2 q_2\rVert_0\leqslant \lVert B_{l_j}\rVert_{0}^{2}\epsilon^{\frac{1}{2}}\times (\lvert c_j\rvert+\epsilon_{l_j}^{\frac{1}{4}})\leqslant C_0\epsilon^{\frac{1}{2}}.$$
Therefore, if we pick $$\epsilon \geqslant C_0 l_j^{-\frac{k}{16}},$$ then $$\lVert Z \rVert_0 \leqslant 1+C_0\epsilon^{\frac{1}{2}}.$$

Now we can prove our main theorem. Notice that it is enough to consider the case when $\epsilon=\lVert \mathcal{B}-\mathcal{A}\rVert_0$ is sufficiently small. Then $\epsilon$ satisfies
\begin{equation}\label{er}
I_j: C_0 l_j^{-\frac{k}{16}}\leqslant \epsilon \leqslant c{l_j^{-\frac{k}{40}}}
\end{equation}
for some choice of $j$ and the following inequality is fulfilled:
$$
C_0 l_j^{-\frac{k}{16}}\leqslant c{l_j^{-\frac{k}{20}}} = c{l_{j+1}^{-\frac{k}{40}}}.
$$
This is essential as we are able to cover all the small $\epsilon$ tending to zero by the interval $I_j$ above.

Let $$\tilde{\mathcal{B}}(\theta)=W(\theta+\alpha)\mathcal{B}(\theta)W(\theta)^{-1},$$ then
$$
\lVert \tilde{\mathcal{B}}\rVert_0 \leqslant \lVert Z \rVert_0 +\lVert W\rVert_0 ^2\lVert \mathcal{B}-\mathcal{A}\rVert_0\leqslant 1+C_0\epsilon^{\frac{1}{2}}.
$$
Thus $$L(\alpha,\mathcal{B})=L(\alpha,\tilde{\mathcal{B}})\leqslant \ln\lVert \tilde{\mathcal{B}}\rVert_0\leqslant C_0\epsilon^{\frac{1}{2}}.$$
Respectively, $$L(\alpha,\mathcal{A})=L(\alpha,Z)\leqslant \ln \lVert Z \rVert_0 \leqslant C_0\epsilon^{\frac{1}{2}}.$$ It implies that $$\lvert L(\alpha,\mathcal{A})-L(\alpha,\mathcal{B})\rvert \leqslant C_0\lVert \mathcal{B}-\mathcal{A}\rVert_0^{\frac{1}{2}}.$$
where $C_0$ depends on $d,\kappa,\tau.$ This finishes the proof of Theorem $\ref{thm1}$.
\end{pf}

\subsection{Proof of Theorem 1.2}
\begin{pf}
We rewrite the Schr\"{o}dinger cocycle as
$$
S_E^{\lambda V}=A+F(\theta),
$$
where
$$
A=\begin{pmatrix} -E & -1 \\ 1 & 0 \end{pmatrix}, \ \ F(\theta)=\begin{pmatrix} \lambda V(\theta) & 0 \\ 0 & 0 \end{pmatrix}.
$$
Notice that we only need to consider $$E\in\Sigma_{\alpha,\lambda V}\subset[-2+\inf \lambda V,2+\sup \lambda V],$$ otherwise $S_E^{\lambda V}$ is uniformly hyperbolic and $L(E)$ (thus $N_{\lambda V,\alpha}(E)$) is automatically $\frac{1}{2}$-H\"older continuous. So by the assumption on $\lambda$, we have $\lVert A\rVert\leqslant 3$ and so does $A^{-1}$.
Clearly, the smallness condition on $\lambda$ is namely the smallness condition on $F$. Thus if we write that
$$
A+F(\theta)=Ae^{f(\theta)},
$$
the assumptions of Theorem $\ref{thm3}$ are naturally fulfilled as
$$
\lVert f(\theta) \rVert_k\leqslant \lVert A^{-1}\rVert \lVert F(\theta)\rVert_k.
$$
Therefore $C^{k, k_0}$ almost reducibility of $(\alpha,S_E^{\lambda V})$ is a straight result of Theorem $\ref{thm3}$.
As a corollary of  Theorem $\ref{thm1}$, we obtain that
$$
L(E+i\epsilon)-L(E)\leqslant C_0\epsilon^{\frac{1}{2}},\ \ E\in \Sigma,
$$
where $C_0$ does not depend on $E$. 

By the Thouless formula $$L(E)=\int \ln\lvert E-E^{'}\rvert dN_{\lambda V,\alpha}(E^{'}),$$
we have
\begin{flalign*}
L(E+i\epsilon)-L(E)=&\frac{1}{2}\int \ln (1+\frac{\epsilon^2}{(E-E^{'})^2})dN_{\lambda V,\alpha}(E^{'})\\
\geqslant &c'(N_{\lambda V,\alpha}(E+\epsilon)-N_{\lambda V,\alpha}(E-\epsilon))
\end{flalign*}
for every $\epsilon>0$ and $c'$ is a numerical constant.
Thus $$N_{\lambda V,\alpha}(E+\epsilon)-N_{\lambda V,\alpha}(E-\epsilon)\leqslant \frac{C_0}{c'}\epsilon^{\frac{1}{2}}$$ for every $\epsilon>0,E\in \Sigma$.
Notice that $N_{\lambda V,\alpha}$ is locally constant in the gaps, this means precisely that $N_{\lambda V,\alpha}$ is $\frac{1}{2}$-H\"older continuous. This finishes the proof of Theorem $\ref{cor1}$.
\end{pf}

\section{Appendix: Proof of Lemma \ref{lem2}}

In this appendix, we prove Lemma \ref{lem2} by the following quantitative Implicit Function Theorem:
\begin{Theorem}\label{0} \cite{BerBia,deimling} Let $X,Y,Z$ be Banach spaces, $U\subset X$ and $V\subset Y$ neighborhoods of $x_0$ and $y_0$ respectively. Fix $s,\delta>0$ and define $B_s(x_0)=\{x\in X\mid \lVert x-x_0\rVert\leqslant s \}, B_{\delta}(y_0)=\{y\in Y\mid \lVert y-y_0\rVert\leqslant \delta \}.$ Let $\Psi\in C^1(U\times V,Z)$ and $B_s(x_0)\times B_{\delta}(y_0)\subset U\times V$. Suppose also that $\Psi(x_0,y_0)=0$, and that $D_y\Psi(x_0,y_0)\in \mathcal{L}(Y,Z)$ is invertible. If
\begin{equation}\label{1}
\sup_{\substack{\overline{B_s(x_0)}}}\lVert \Psi(x,y_0)\rVert_Z\leqslant \frac{\delta}{2\lVert ({D_y\Psi(x_0,y_0)}^{-1})\rVert},
\end{equation}
\begin{equation}\label{2}
\sup_{\substack{\overline{B_s(x_0)}\times \overline{B_{\delta}(y_0)}}}\lVert Id_Y-(D_y\Psi(x_0,y_0))^{-1}D_y\Psi(x,y) \rVert_{\mathcal{L}(Y,Y)}\leqslant \frac{1}{2},
\end{equation}
then there exists $y\in C^1(B_s(x_0),\overline{B_{\delta}(y_0)})$ such that $\Psi(x,y(x))=0$.
\end{Theorem}

With Theorem \ref{0} in hand,  now we can prove Lemma \ref{lem2} easily.  We construct the nonlinear functional
$$
\Psi:\mathcal{B}_r^{nre}(\eta) \times C^{\omega}_r(\T^d,su(1,1))\rightarrow \mathcal{B}_r^{nre}(\eta)
$$
by
$$
\Psi(Y,g)=\mathbb{P}_{nre} \ln(e^{A^{-1}Y(\theta+\alpha)A}e^{g(\theta)}e^{-Y(\theta)})
$$
An immediate check reveals that $$\Psi(0,0)=0, \ \ \lVert \Psi(0,g)\rVert\leqslant \lvert g\rvert_r.
$$
and
\begin{flalign*}
&\Psi(Y+Y',g)-\Psi(Y,g)\\
=&\mathbb{P}_{nre} \ln(e^{A^{-1}(Y(\theta+\alpha)+Y'(\theta+\alpha))A}e^{g(\theta)}e^{-(Y(\theta)+Y'(\theta))})\\
&-\mathbb{P}_{nre} \ln(e^{A^{-1}Y(\theta+\alpha)A}e^{g(\theta)}e^{-Y(\theta)})\\
=&\mathbb{P}_{nre} \ln(e^{A^{-1}(Y(\theta+\alpha)+Y'(\theta+\alpha))A}e^{g(\theta)}e^{-(Y(\theta)+Y'(\theta))})\\
&-\mathbb{P}_{nre}\ln(e^{A^{-1}(Y(\theta+\alpha)+Y'(\theta+\alpha))A}e^{g(\theta)}e^{-Y(\theta)})\\
&+\mathbb{P}_{nre}\ln(e^{A^{-1}(Y(\theta+\alpha)+Y'(\theta+\alpha))A}e^{g(\theta)}e^{-Y(\theta)})\\
&-\mathbb{P}_{nre} \ln(e^{A^{-1}Y(\theta+\alpha)A}e^{g(\theta)}e^{-Y(\theta)})
\end{flalign*}
To make further computations, we need the fact that if $A,B,C$ are small $sl(2,\R)$ matrices, then there exists $D,E\in sl(2,\R)$ such that
$$
e^{A}e^{B}e^{C}=e^{D+E}
$$
and
$$
D=A+B+C
$$
where $E$ is a sum of terms of order at least 2 in $A,B,C$.
Also, the famous Baker-Campbell-Hausdorff Formula shows that
\begin{equation}
\ln(e^X e^Y)=X+Y+\frac{1}{2}[X,Y]+\frac{1}{12}([X,[X,Y]+[Y,[Y,X]])+\cdots,
\end{equation}
where $[X,Y]=XY-YX$ denotes the Lie Bracket and $\cdots$ denotes the sum of higher order terms.

Therefore, we can compute that
\begin{flalign*}
&\mathbb{P}_{nre} \ln(e^{A^{-1}(Y(\theta+\alpha)+Y'(\theta+\alpha))A}e^{g(\theta)}e^{-(Y(\theta)+Y'(\theta))})\\
-&\mathbb{P}_{nre}\ln(e^{A^{-1}(Y(\theta+\alpha)+Y'(\theta+\alpha))A}e^{g(\theta)}e^{-Y(\theta)})\\
=&\mathbb{P}_{nre} \ln(e^{A^{-1}(Y(\theta+\alpha)+Y'(\theta+\alpha))A}e^{g(\theta)}e^{-Y(\theta)}e^{-Y''(\theta)})\\
-&\mathbb{P}_{nre}\ln(e^{A^{-1}(Y(\theta+\alpha)+Y'(\theta+\alpha))A}e^{g(\theta)}e^{-Y(\theta)})\\
=&\mathbb{P}_{nre} \ln(e^{D+E}e^{-Y''(\theta)})-\mathbb{P}_{nre}\ln(e^{D+E})\\
=&\mathbb{P}_{nre}(D+E-Y''+\frac{1}{2}[D+E,-Y'']+\cdots)-\mathbb{P}_{nre}(D+E)\\
=&\mathbb{P}_{nre}(-Y''+\frac{1}{2}[D+E,-Y'']+\cdots)
\end{flalign*}
where
$$
Y''(\theta)=Y'(\theta)+\mathcal{O}(Y(\theta))Y'(\theta),
$$
$$
D(\theta)=A^{-1}(Y(\theta+\alpha)+Y'(\theta+\alpha))A+g(\theta)-Y(\theta),
$$
and $E$ is a sum of terms of order at least 2 in $A^{-1}(Y(\theta+\alpha)+Y'(\theta+\alpha))A,g(\theta),-Y(\theta)$.

Similarly, we have
\begin{flalign*}
&\mathbb{P}_{nre}\ln(e^{A^{-1}(Y(\theta+\alpha)+Y'(\theta+\alpha))A}e^{g(\theta)}e^{-Y(\theta)})\\
-&\mathbb{P}_{nre} \ln(e^{A^{-1}Y(\theta+\alpha)A}e^{g(\theta)}e^{-Y(\theta)})\\
=&\mathbb{P}_{nre}\ln(e^{Y'''}e^{A^{-1}Y(\theta+\alpha)A}e^{g(\theta)}e^{-Y(\theta)})\\
-&\mathbb{P}_{nre} \ln(e^{A^{-1}Y(\theta+\alpha)A}e^{g(\theta)}e^{-Y(\theta)})\\
=&\mathbb{P}_{nre} \ln(e^{Y'''}e^{F+H})-\mathbb{P}_{nre}\ln(e^{F+H})\\
=&\mathbb{P}_{nre}(Y'''+F+H+\frac{1}{2}[Y''',F+H]+\cdots)-\mathbb{P}_{nre}(F+H)\\
=&\mathbb{P}_{nre}(Y'''+\frac{1}{2}[Y''',F+H]+\cdots)
\end{flalign*}
where
$$
Y'''(\theta+\alpha)=A^{-1}Y'(\theta+\alpha)A+\mathcal{O}(A^{-1}Y(\theta+\alpha)A)\times A^{-1}Y'(\theta+\alpha)A,
$$
$$
F(\theta)=A^{-1}Y(\theta+\alpha)A+g(\theta)-Y(\theta)
$$
and $H$ is a sum of terms of order at least 2 in $A^{-1}(Y(\theta+\alpha)A,g(\theta),-Y(\theta)$.

By the definition of Frechet differential, we only need to consider the linear terms of $\Psi(Y+Y',g)-\Psi(Y,g)$, thus we have
\begin{flalign*}
D_Y\Psi(Y,g)(Y')&=\mathbb{P}_{nre}(A^{-1}Y'(\theta+\alpha)A+\mathcal{O}(A^{-1}Y(\theta+\alpha)A)\times A^{-1}Y'(\theta+\alpha)A\\
&+\frac{1}{2}[Y''',F+H]+\cdots)\\
&+\mathbb{P}_{nre}(-Y'(\theta)-\mathcal{O}(Y(\theta))Y'(\theta)+\frac{1}{2}[F+H',-Y'']+\cdots),
\end{flalign*}
where $H'$ is a sum of terms of order at least 2 in $A^{-1}(Y(\theta+\alpha)A,g(\theta),-Y(\theta)$. Moreover, the first \textquotedblleft $\cdots$\textquotedblright denotes the sum of terms which have order at least $2$ in $F+H$ but only order $1$ in $Y'''$. The second \textquotedblleft $\cdots$\textquotedblright denotes the sum of terms which have order at least $2$ in $F+H'$ but only order $1$ in $Y''$.

Let $Y=0$ and $g=0$, then all the Lie brackets vanish. So we immediately obtain
\begin{flalign*}
D_Y\Psi(0,0)(Y')&=\mathbb{P}_{nre}(A^{-1}Y'(\theta+\alpha)A)+\mathbb{P}_{nre}(-Y'(\theta))\\
&=A^{-1}Y'(\theta+\alpha)A-Y'(\theta)
\end{flalign*}
Thus
\begin{flalign*}
\lVert D_Y\Psi(0,0)(Y')\rVert &\geqslant \lvert A^{-1}Y'(\theta+\alpha)A-Y'(\theta)\rvert_r \\
&\geqslant \eta\lvert Y'\rvert_r.\\
\end{flalign*}
So we have
$$
\lVert (D_Y\Psi(0,0))^{-1}\rVert\leqslant \eta^{-1}.
$$
For our purpose, we set $s=\epsilon,\delta=\epsilon^\frac{1}{2}$ and $\eta \geqslant 13\lVert A\rVert^2\epsilon^{\frac{1}{2}}$. Then we have
$$
2\times\sup_{\substack{\overline{B_{s}(0)}}}\lVert \Psi(0,g)\rVert \times \lVert (D_Y\Psi(0,0))^{-1}\rVert\leqslant 2\times \epsilon\times \frac{1}{13\lVert A\rVert^2}\epsilon^{-\frac{1}{2}}\leqslant \epsilon^{\frac{1}{2}}=\delta,
$$
then $(\ref{1})$ is fulfilled.

On the other hand, direct computation shows that
\begin{flalign*}
&D_Y\Psi(Y,g)(Y')-D_Y\Psi(0,0)(Y')\\
=&\mathbb{P}_{nre}(\mathcal{O}(A^{-1}Y(\theta+\alpha)A)\times A^{-1}Y'(\theta+\alpha)A+\frac{1}{2}[Y''',F+H]+\cdots)\\
+&\mathbb{P}_{nre}(-\mathcal{O}(Y(\theta))Y'(\theta)+\frac{1}{2}[F+H',-Y'']+\cdots).
\end{flalign*}
Therefore, we have
\begin{flalign*}
\sup_{\substack{\overline{B_s(0)}\times \overline{B_{\delta}(0)}}}\lVert D_Y\Psi(Y,g)(Y')-D_Y\Psi(0,0)(Y')\rVert & \leqslant 6(\lVert A\rVert^2\lvert Y\rvert_r+\lVert A\rVert^2\lvert g\rvert_r)\lvert Y'\rvert\\
&\leqslant  6(\lVert A\rVert^2\delta+\lVert A\rVert^2s)\lvert Y'\rvert\\
&\leqslant  6(\lVert A\rVert^2 \epsilon^{\frac{1}{2}}+\lVert A\rVert^2\epsilon)\lvert Y'\rvert,
\end{flalign*}
which implies
$$
\sup_{\substack{\overline{B_s(0)}\times \overline{B_{\delta}(0)}}}\lVert D_Y\Psi(0,0)-D_Y\Psi(Y,g)\rVert\leqslant 6(\lVert A\rVert^2 \epsilon^{\frac{1}{2}}+\lVert A\rVert^2\epsilon).
$$
Thus we have
\begin{flalign*}
&\sup_{\substack{\overline{B_s(0)}\times \overline{B_{\delta}(0)}}}\lVert Id_{\mathcal{B}_r^{nre}(\eta)}-(D_Y\Psi(0,0))^{-1}\times D_Y\Psi(Y,g) \rVert\\
\leqslant & \sup_{\substack{\overline{B_s(0)}\times \overline{B_{\delta}(0)}}}\lVert D_Y\Psi(0,0)-D_Y\Psi(Y,g)\rVert \times \lVert (D_Y\Psi(0,0))^{-1}\rVert\\
\leqslant & 6(\lVert A\rVert^2 \epsilon^{\frac{1}{2}}+\lVert A\rVert^2\epsilon) \times \frac{1}{13\lVert A\rVert^2}\epsilon^{-\frac{1}{2}}\\
\leqslant & \frac{1}{2},
\end{flalign*}
which satisfies $(\ref{2})$.
By Theorem $\ref{0}$, for $\lvert g\rvert_r\leqslant \epsilon$ and $\eta \geqslant 13\lVert A\rVert^2\epsilon^{\frac{1}{2}}$, there exists $\lvert Y\rvert_r\leqslant \epsilon^{\frac{1}{2}}$ such that $\Psi(Y,g)=0$, i.e.
$$
e^{A^{-1}Y(\theta+\alpha)A}e^{g(\theta)}e^{Y(\theta)}=e^{g^{re}(\theta)},
$$
which is equivalent to 
$$
e^{Y(\theta+\alpha)}Ae^{g(\theta)}e^{Y(\theta)}=Ae^{g^{re}(\theta)}
$$
and it is easy to check $\lvert g^{re} (\theta)\rvert_r\leqslant 2\epsilon$. This finishes the proof of Lemma \ref{lem2}.

\section{Acknowledgements}
C. Chavaudret was supported by  the ANR "BEKAM" and  the ANR "Dynamics and CR Geometry". J. You  was partially supported by NSFC grant (11471155) and
973 projects of China (2014CB340701).   Q. Zhou was partially supported by \textquotedblleft Deng Feng Scholar Program B\textquotedblright of Nanjing University, Specially-appointed professor programe of Jiangsu province and NSFC grant (11671192).

\end{document}